\documentclass{amsart}
\usepackage{poisson}

\title{A Poisson relation for conic manifolds}
\author{Jared Wunsch}
\date{\today}
\address{Department of Mathematics\\Northwestern University\\}
\email{jwunsch@math.sunysb.edu}

\begin{document}

\begin{abstract}
Let $X$ be a compact Riemannian manifold with conic singularities, i.e.\ a
Riemannian manifold whose metric has a conic degeneracy at the
boundary. Let $\Lap$ be the Friedrichs extension of the Laplace-Beltrami
operator on $X.$ There are two natural ways to define geodesics passing
through the boundary: as ``diffractive'' geodesics which may emanate from
$\pa X$ in any direction, or as ``geometric'' geodesics which must enter
and leave $\pa X$ at points which are connected by a geodesic of length
$\pi$ in $\pa X.$ Let $\dif=\{0\} \cup \{\pm \text{lengths of closed
diffractive geodesics}\}$ and $\geo=\{0\} \cup \{\pm \text{lengths of }$
$\text{closed geometric geodesics}\}.$ We show that
$$
\Tr \cos t \sqrt\Lap \in \cC^{-n-0}(\RR) \cap \cC^{-1-0}(\RR\backslash
\geo) \cap \CI(\RR\backslash \dif).$$ This generalizes a classical result
of Chazarain and Duistermaat-Guillemin on boundaryless manifolds, which in
turn follows from Poisson summation in the case $X=S^1.$
\end{abstract}

\maketitle

\section{Introduction}
Let $(X,g)$ be a conic manifold, that is to say, a manifold $X$ with
boundary, endowed with a Riemannian metric that takes the form
$$
g=dx^2+x^2 h
$$
in a neighborhood of $\pa X$, with $x$ a boundary defining function and $$h
\in \CI(X; \Sym^2 T^*X)$$ nondegenerate in the sense that $h_0 \equiv h
\restrictedto_{x=0}$ is a metric.  Trivial examples of conic metrics are
obtained by blowing up isolated points on a Riemannian manifold.  Less
trivial examples include the \emph{product cones} $\overline{\RR_+}\times
Y$ endowed with metric $dx^2+x^2 h(y,dy)$ with $h$ a metric on $Y$.  A
conic manifold in the sense defined here can be regarded as the
desingularization, or blow-up, of a manifold with conic singularities.  The
degeneracy of the metric at the boundary makes each boundary component look
geometrically like a separate cone point.  One geometric consequence of the
degeneracy of the metric is the most geodesics fail to strike $\pa X$;
indeed through any point close to $\pa X$ there will be precisely one
direction in which the exponential map will reach the boundary in short
time (see \cite{Melrose-Wunsch1}).

Let $\Lap$ be the Friedrichs extension of the (nonnegative)
Laplace-Beltrami operator on $X$.  The propagation of singularities for
solutions of the wave equation 
$$
(D_t^2-\Lap) u =0
$$ on conic manifolds was studied by the author and R.~Melrose in
\cite{Melrose-Wunsch1}.  The results of that paper showed that while
singularities propagate along geodesics over $X^\circ$ as described by
H\"ormander's theorem \cite{MR48:9458}, a singularity striking one
component of $\pa X$ will, in general, give rise to a ``diffracted wave''
emanating spherically from that component of the boundary.  The singularity
of the diffracted wave, measured away from points which are ``geometrically
related'' to the incoming singularity by time-$\pi$ geodesic flow in $\pa
X$, will in many cases be weaker than the incident singularity; for this
difference in regularity to manifest itself, it suffices that the incident
singularity be not too directly \emph{focused} upon the boundary, in a
sense made precise in \cite{Melrose-Wunsch1}.  A plausible explanation for
this distinction is as follows: If we take a sequence of geodesics in the
interior of $X$ passing arbitrarily close, in the limit, to $\pa X$, then
there exists a subsequence approaching a continuous path consisting of: a
geodesic striking the boundary; followed by a geodesic \emph{within} the
boundary of length $\pi$; followed by a geodesic emanating from the
boundary.  This outgoing geodesic is thus geometrically related in the
sense described above to the incoming one.  Thus strong singularities
propagate through the boundary as one would expect from their limiting
behavior in the interior; the (often weaker) diffractive wave, by contrast,
sends singularities out from the boundary along rays which are
geometrically inacessible and lie in this sense in a shadow region.

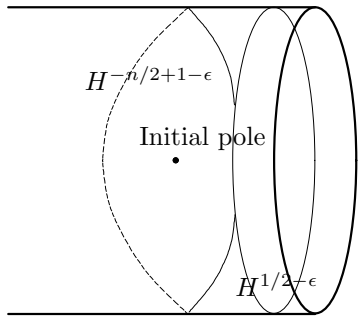
\begin{figure}
\setlength{\unitlength}{0.00041667in}
\begingroup\makeatletter\ifx\SetFigFont\undefined%
\gdef\SetFigFont#1#2#3#4#5{%
  \reset@font\fontsize{#1}{#2pt}%
  \fontfamily{#3}\fontseries{#4}\fontshape{#5}%
  \selectfont}%
\fi\endgroup%
{
\begin{picture}(4445,3942)(0,-10)
\put(2142,1965){\blacken\ellipse{68}{68}}
\put(2142,1965){\ellipse{68}{68}}
\thinlines
\put(3387,1965){\ellipse{1050}{3900}}
\thicklines
\put(3912,1965){\ellipse{1050}{3900}}
\path(3912,15)(12,15)
\path(3912,3915)(12,3915)
\thinlines
\path(2888,2670)(2888,2671)(2888,2674)
	(2887,2682)(2886,2697)(2884,2717)
	(2881,2743)(2878,2773)(2875,2805)
	(2871,2837)(2867,2868)(2863,2898)
	(2859,2926)(2855,2951)(2850,2976)
	(2846,2999)(2841,3021)(2835,3043)
	(2829,3065)(2823,3083)(2818,3102)
	(2811,3121)(2804,3141)(2797,3161)
	(2789,3182)(2781,3203)(2772,3225)
	(2763,3247)(2753,3270)(2742,3292)
	(2732,3315)(2721,3337)(2710,3360)
	(2699,3381)(2687,3403)(2676,3424)
	(2665,3444)(2653,3463)(2642,3482)
	(2630,3501)(2619,3519)(2607,3537)
	(2595,3554)(2582,3572)(2569,3591)
	(2555,3609)(2540,3629)(2524,3650)
	(2506,3671)(2487,3694)(2467,3719)
	(2445,3744)(2423,3770)(2401,3796)
	(2380,3820)(2360,3843)(2342,3864)
	(2327,3880)(2316,3893)(2309,3901)
	(2305,3906)(2303,3908)
\path(2895,1275)(2895,1274)(2895,1271)
	(2894,1263)(2893,1248)(2891,1228)
	(2888,1202)(2885,1172)(2882,1140)
	(2878,1108)(2874,1077)(2870,1047)
	(2866,1019)(2862,994)(2857,969)
	(2853,946)(2848,924)(2842,902)
	(2836,880)(2830,862)(2825,843)
	(2818,824)(2811,804)(2804,784)
	(2796,763)(2788,742)(2779,720)
	(2770,698)(2760,675)(2749,653)
	(2739,630)(2728,608)(2717,585)
	(2706,564)(2694,542)(2683,521)
	(2672,501)(2660,482)(2649,463)
	(2637,444)(2626,426)(2614,408)
	(2602,391)(2589,373)(2576,354)
	(2562,336)(2547,316)(2531,295)
	(2513,274)(2494,251)(2474,226)
	(2452,201)(2430,175)(2408,149)
	(2387,125)(2367,102)(2349,81)
	(2334,65)(2323,52)(2316,44)
	(2312,39)(2310,37)
\path(2292,3904)(2290,3903)(2287,3899)
\path(2280,3893)(2270,3884)(2256,3871)
\path(2238,3854)(2217,3834)(2192,3811)
\path(2165,3785)(2136,3758)(2106,3729)
\path(2076,3700)(2046,3671)(2017,3642)
\path(1988,3614)(1961,3587)(1936,3561)
\path(1911,3536)(1888,3512)(1866,3488)
\path(1846,3465)(1826,3443)(1807,3421)
\path(1789,3399)(1771,3377)(1753,3355)
\path(1736,3332)(1718,3309)(1701,3286)
\path(1684,3261)(1666,3237)(1649,3211)
\path(1631,3185)(1614,3158)(1596,3130)
\path(1578,3102)(1560,3073)(1543,3043)
\path(1525,3013)(1508,2983)(1491,2953)
\path(1474,2922)(1458,2892)(1442,2862)
\path(1427,2832)(1412,2803)(1398,2774)
\path(1385,2746)(1373,2719)(1361,2692)
\path(1350,2666)(1339,2641)(1329,2616)
\path(1320,2592)(1312,2568)(1302,2541)
\path(1294,2515)(1286,2488)(1278,2462)
\path(1272,2434)(1265,2407)(1259,2378)
\path(1254,2348)(1249,2316)(1244,2283)
\path(1239,2249)(1235,2213)(1231,2177)
\path(1227,2141)(1224,2106)(1221,2073)
\path(1218,2044)(1216,2018)(1214,1998)
\path(1213,1983)(1213,1973)(1212,1967)(1212,1965)
\path(2300,12)(2298,13)(2295,17)
\path(2288,23)(2278,32)(2264,45)
\path(2246,62)(2225,82)(2200,105)
\path(2173,131)(2144,158)(2114,187)
\path(2084,216)(2054,245)(2025,274)
\path(1996,302)(1969,329)(1944,355)
\path(1919,380)(1896,404)(1874,428)
\path(1854,451)(1834,473)(1815,495)
\path(1797,517)(1779,539)(1761,561)
\path(1744,584)(1726,607)(1709,630)
\path(1692,655)(1674,679)(1657,705)
\path(1639,731)(1622,758)(1604,786)
\path(1586,815)(1568,844)(1550,873)
\path(1533,903)(1516,933)(1498,964)
\path(1482,994)(1466,1024)(1450,1055)
\path(1435,1084)(1420,1114)(1406,1143)
\path(1393,1171)(1380,1198)(1368,1225)
\path(1357,1251)(1346,1277)(1336,1302)
\path(1327,1326)(1318,1350)(1309,1377)
\path(1300,1404)(1292,1431)(1284,1458)
\path(1277,1486)(1271,1514)(1265,1543)
\path(1259,1574)(1253,1606)(1248,1640)
\path(1243,1675)(1238,1711)(1234,1748)
\path(1229,1785)(1226,1820)(1222,1854)
\path(1219,1884)(1217,1910)(1215,1931)
\path(1214,1947)(1213,1957)(1212,1962)(1212,1965)
\put(1677,2160){Initial pole}
\put(987,2865){$H^{-n/2+1-\epsilon}$}
\put(2900,240){$H^{1/2-\epsilon}$}
\end{picture}
}
\caption{The wavefronts of the fundamental solution $\sin
t\sqrt\Lap/\sqrt\Lap$}\label{figure}
\end{figure}

In the crucial example of the fundamental solution $\sin t \sqrt
\Lap/\sqrt\Lap$ with initial pole chosen sufficiently close to the
boundary, the diffracted wave is $(n-1)/2-\ep$ derivatives smoother than
the main singularity (for all $\ep>0$), and is furthermore conormal to the
radial surface emanating from $\pa X$.  Figure~\ref{figure} shows the main
front and the diffractive wave for such a solution, labeled by their
Sobolev regularity.  In this case the diffractive region is the part of the
spherical wave emanating from the boundary that does not lie in the closure
of the main front.

In the case of product cones, these propagation results follow from the
explicit construction of the fundamental solution carried out by Cheeger
and Taylor \cite{Cheeger-Taylor1,Cheeger-Taylor2} See also the work of
Gerard-Lebeau \cite{Gerard-Lebeau2} for an analogous result in analytic
category in the setting of manifolds with boundary.

Let $X$ be a conic manifold and let $Y_i$ denote the boundary components of
$X$, with $i=1,\dots,K$.  Based on the two different kinds of propagation
of singularities described above, we now make two different definitions of
geodesics passing through $\pa X.$
\begin{definition}
A \emph{diffractive geodesic} on $X$ is a union of a finite number of
closed, oriented geodesic segments $\gamma_1, \dots, \gamma_N$ in $X$
such that all end points except possibly the initial point in $\gamma_1$ and the
final point of $\gamma_N$ lie in $\pa X$, and $\gamma_i$ ends at the same
boundary component at which $\gamma_{i+1}$ begins, for $i=1, \dots, N-1$.

A \emph{geometric geodesic} is a diffractive geodesic such that in
addition, the final point of $\gamma_i$ and the initial point of
$\gamma_{i+1}$ are connected by a geodesic of length $\pi$ in $\pa X$
(w.r.t.\ the metric $h_0=h\restrictedto_{\pa X}$) for $i=1,\dots, N-1$.

A diffractive geodesic is \emph{closed} if the initial point of $\gamma_1$
and the final point of $\gamma_N$ coincide in $X^\circ$ and if $\gamma_1$,
$\gamma_N$ have the same tangent there, or if the initial point of
$\gamma_1$ and the final point of $\gamma_N$ lie in the same component of
$\pa X.$

A geometric geodesic is \emph{closed} if the initial point of $\gamma_1$
and the final point of $\gamma_N$
coincide in $X^\circ$ and if $\gamma_1$, $\gamma_N$ have the
same tangent there, or if the initial point of $\gamma_1$ and the final
point of $\gamma_N$ lie in $\pa X$ and are connected by a geodesic of
length $\pi$ in $\pa X.$

The \emph{length} of a diffractive geodesic is the sum of the lengths of
its segments $\gamma_i$.

Let $$\dif= \{\pm \text{ lengths of closed diffractive
geodesics}\} \cup \{0\}$$
and
$$\geo= \{\pm \text{ lengths of closed geometric
geodesics}\} \cup \{0\}.$$
\end{definition}

Let $k \in \NN$.  Let $u$ be a distribution on $\RR$.  We write $u \in
\cC^{-k-0}(\RR)$ if $D_t^{-k-\delta} \phi u \in L^\infty_{\loc}(\RR)$ for
all $\delta>0$ and $\phi \in \mathcal{C}_c^\infty(\RR).$ Note that if $u$
has compact support, this coincides with $u \in C_*^{-k-\delta}(\RR)$ for
all $\delta>0$ in the notation of \cite{Taylor7}, Section 13.8.

This paper is devoted to establishing the following ``Poisson relation''
(so-called because it generalized part of the Poisson summation formula
from the case in which $X$ is the boundaryless manifold $S^1$):
\begin{thm}
Let $\Lap$ be the Friedrichs extension of the Laplacian on a compact conic
manifold.  Then
$$
\Tr \cos t\sqrt\Lap \in \cC^{-n-0}(\RR) \cap \cC^{-1-0}(\RR\backslash
\geo) \cap \CI(\RR\backslash \dif).
$$
\end{thm}
This generalizes the classical result on a boundaryless manifold, due to
Chazarain \cite{Chazarain1} and Duistermaat-Guillemin
\cite{Duistermaat-Guillemin1}, that the trace is smooth away from
$\{0\}\cup \{\pm \text{lengths of closed geodesics}\}$.  Note that when the
closed geodesics on a compact, boundaryless manifold are isolated and
nondegenerate, the results of Chazarain and Duistermaat-Guillemin show that
$\Tr \cos t \sqrt \Lap \in \cC^{-1-0}(\RR\backslash \{0\})$.  It thus seems
possible that subject to an appropriate nondegeneracy assumption,
regularity in $\cC^{-1-0}$ might also hold at lengths of all closed
geometric geodesics.  Thus it is possible that the strength of the
singularities in the wave trace allows us to distinguish lengths of
geometric geodesics from diffractive ones only when the former are actually
\emph{degenerate}.

It is natural to ask whether a full trace formula of the type proved by
Chazarain and Duistermaat-Guillemin, with asymptotic expansions at the
singularities, holds on conic manifolds.  Unfortunately the tools developed
in \cite{Melrose-Wunsch1} are insufficiently constructive to demonstrate
such an expansion except perhaps at $t=0$.  A more accessible future
direction would be to extend the results of this paper to give a filtration
of the singularities according to how many diffractive interactions with
the boundary a closed geodesic undergoes: subject to some nondegeneracy
assumptions, more diffractions should result in a weaker singularity.

The author is grateful to Richard Melrose for helpful discussions and to an
anonymous referee for substantial improvements to the exposition.  This
research was supported in part by NSF grant DMS-0100501.

\section{Geodesic flow}\label{section:geodesic}
We begin by describing the lifts of the diffracted and geometric geodesics
to (a rescaled version of) the cotangent bundle.

In \cite{Melrose-Wunsch1}, it is shown that by appropriate choice of
product structure $(x,y)$ near a boundary component $Y_i$, we may bring any
conic metric to the reduced form
\begin{equation}
g=dx^2+x^2 h(x,y,dy),
\label{semiproduct}
\end{equation}
i.e.\ we can reduce $h$ to a smooth family in $x$ of metrics on $Y_i.$ The
existence of this normal form is equivalent to the existence of a fibration
of a neighborhood of $Y_i$ by short geodesics reaching $Y_i$: with the
metric in the form \eqref{semiproduct} the geodesics reaching $Y_i$ are
precisely those of the form $y=y_0,\ x=x_0-t.$  \emph{Henceforth, we work in such
product coordinates.}

Let $\Tbstar X$ denote the b-cotangent bundle of $X$, i.e.\ the dual of the
bundle whose sections are smooth vector fields tangent to $\pa X$.  Let
$\Sbstar X$ denote the corresponding sphere bundle.  Let $\xi \, dx/x+
\eta \cdot dy$ denote the canonical one-form on $\Tbstar X$.

Let $K_g$ be the Hamilton vector field for
$g/2=(\xi^2+h(x,y,\eta))/(2x^2)$, the symbol of $\Lap/2$ on $\Tbstar X;$
note that $K_g$ is merely the geodesic spray in $\Tbstar X$ with velocity
$\sqrt g.$ It is convenient to rescale this vector field so that it is both
tangent to the boundary of $X$ and homogeneous of degree zero in the
fibers. Near a boundary
component $Y_i$, for a metric in the reduced form \eqref{semiproduct}, we
have (see \cite{Melrose-Wunsch1})
\begin{equation}
K_g=x^{-2} \bigg( H_{Y_i}(x) + \big(\xi^2+h(x,y,\eta)+ \frac x2 \frac{\pa
h}{\pa x}\big) \pa_\xi +\xi x \pa _x \bigg),
\label{spray}
\end{equation}
where $H_{Y_i}(x)$ is the geodesic spray in $Y_i$ with respect to the
family of metrics $h(x, \cdot)$.  Hence the desired rescaling is
$$
Z=\frac{x}{\sqrt g} K_g.
$$ By the homogeneity of $Z,$ if we radially compactify the fibers of the
cotangent bundle and identify $\Sbstar X$ with the ``sphere at infinity''
then $Z$ is tangent to $\Sbstar X,$ and may be restricted to it.
Henceforth, then, we let $Z$ denote the restriction of $(x/\sqrt g) K_g$ to
the \emph{compact} manifold $\Sbstar X$ on which the coordinates $\xi,
\eta$ have been replaced by the (redundant) coordinates $$(\bar\xi,
\bar\eta)=\left(\xi/\sqrt{\xi^2+h(\eta)},
\eta/\sqrt{\xi^2+h(\eta)}\right).$$ $Z$ vanishes only at certain points
$x=\bar\eta=0$ over $\pa X$, hence the closures of maximally extended
integral curves of this vector field can only begin and end over $\pa X$.
Since $Z$ is tangent to the boundary, such integral curves either lie
entirely over $\pa X$ or lie over $\pa X$ only at their limit points.
Interior and boundary integral curves can meet only at limit points in
$\{x=\bar\eta=0\} \subset \Sbstar X.$

It is helpful in studying the integral curves of $Z$ to introduce the
following way of measuring their lengths: Let $\gamma$ be an integral curve
of $Z$ over $X^\circ$.  Let $k$ denote a Riemannian metric on $\Sbstar
X^\circ$ such that $k(Z,Z)=1$.  Let
\begin{equation}
\omega= x k(\cdot, Z) \in \Omega^1(\Sbstar X).
\label{omega}
\end{equation}
Then
$$
\int_\gamma \omega= \int_\gamma x k(d\gamma/ds, Z)\, ds =\int_\gamma \frac x{\sqrt g}
k(K_g, Z)\, ds=\int_\gamma ds =\text{length}(\gamma) 
$$
where $s$ parametrizes $\gamma$ as an integral curve of $K_g/\sqrt{g}$, the
unit speed geodesic flow.  With this motivation in mind, we now define, for
each $t\in\RR_+$, two relations in $\Sbstar X$, a ``geometric'' and a
``diffractive'' relation.  These correspond to the two different
possibilities for geodesic flow through the boundary.
\begin{definition}\label{defn:relations}
Let $p, q \in \Sbstar X$.  We write
$$
p \gtilde q
$$
if there exists a \emph{continuous,} piecewise smooth curve $\gamma: [0,1] \to
\Sbstar X$ with $\gamma(0)=p$, $\gamma(1)=q$, such that $[0,1]$ can be
decomposed into a finite union of closed subintervals $I_j$, intersecting at
their endpoints, where
\begin{enumerate}
\item on each $I_j^\circ$, $\gamma$ is a
(reparametrized) positively oriented integral curve of $Z$ in $\Sbstar X,$
\item
On successive intervals $I_j$ and $I_{j+1}$, interior and boundary curves
alternate,
\item
$\int_\gamma \omega =t$, with $\omega$ as defined
in \eqref{omega}.
\end{enumerate}

We write
$$
p \dtilde q
$$
if there exists a piecewise smooth (not necessarily continuous) curve
$\gamma: [0,1] \to \Sbstar X$ with $\gamma(0)=p$, $\gamma(1)=q$, such that
$[0,1]$ can be decomposed into a finite union of closed subintervals $I_j$,
intersecting at their endpoints, where
\begin{enumerate}
\item on each $I_j^\circ$, $\gamma$ is a
(reparametrized) positively oriented integral curve of $Z$ in $\Sbstar X^\circ,$
\item
the final point of $\gamma$ on $I_j$ and the initial point of $\gamma$ on
$I_{j+1}$ lie over the same component of $\pa X,$
\item
$\int_\gamma
\omega =t$.
\end{enumerate}
\end{definition}

Integral curves of $Z$ over $X^\circ$ are lifts of geodesics in $X^\circ,$ and it
follows from \eqref{spray} that the maximally extended integral curves of
$Z$ in $\Sbstar_{\pa X} X$ are lifts of geodesics of length $\pi$ in $\pa
X$ (see \cite{Melrose-Wunsch1} for details), hence:
\begin{proposition}
$p\gtilde q$ iff $p$ and $q$ are connected by a (lifted) geometric geodesic of
length $t$.

$p\dtilde q$ iff $p$ and $q$ are connected by a (lifted) diffractive geodesic of
length $t$.
\end{proposition}

\begin{proposition}
The sets $\{(p,q,t): p\gtilde q\}$ and $\{(p,q,t): p\dtilde q\}$ are closed
subsets of $\Sbstar X \times \Sbstar X \times \RR_+$.
\label{prop:closed}
\end{proposition}
\begin{proof}
Let $\gamma_i$ be a sequence of ``geometric'' curves as described in the
first part of Definition~\ref{defn:relations} with initial and final points
$p_i, q_i$ where the $p_i$ converge to $p$ and $q_i$ to $q$, and with
$\int_{\gamma_i} \omega \to t$.  Parametrize $\gamma_i=\gamma_i(s)$ with
respect to length in some fixed, nondegenerate Riemannian metric $G$ on
$\Sbstar X$.  The lengths of $\gamma_i$ w.r.t.\ $G$ are uniformly bounded,
since only a finite number of segments are involved (note that the number
of smooth integral curves into which $\gamma_i$ decomposes is a priori
bounded because each interior segment has length bounded below by
$\min_{j,k} d_g(Y_j, Y_k)$).  The $\gamma_i(s)$ are then equicontinuous,
hence applying Ascoli's theorem we may pass to a subsequence that converges
uniformly to a path $\gamma(s)$.  This must be a geometric curve as in
Definition~\ref{defn:relations}, since it is continuous and away from the
singular points of $Z$, a limit of integral curves is an integral curve.
At the set $\{x=\bar \eta=0\} \subset \Sbstar_{\pa X} X$ where $Z$ vanishes, it
vanishes nondegenerately (see \S1 of \cite{Melrose-Wunsch1}), so we may
arrange that the coefficients of the metric $k$ in \eqref{omega} be bounded
by some multiple of $(x^2+\bar\eta^2)^{-1}.$ Since $xZ$ vanishes
quadratically at $\{x= \bar\eta=0\}$, $\omega =x k(Z, \cdot)$ is in
fact in $L^\infty$.  Therefore $\int_\gamma \omega = \lim \int_{\gamma_i}
\omega = t$, as desired.

The result for the diffractive relation follows via a similar argument.
\end{proof}

As a consequence, we have:
\begin{proposition}
For any compact conic manifold, $\dif$ and $\geo$ are closed subsets of
$\RR$.  For any $T\in \RR_+$, the set of lifts of closed diffractive resp.\
geometric geodesics of length $T$ is a compact subset of $\Sbstar X$.
\end{proposition}

\section{Results from \cite{Melrose-Wunsch1}}

In \cite{Melrose-Wunsch1}, the domains of a range of powers of the Friedrichs
Laplacian on a conic manifold $X$ are identified in terms of certain Sobolev
spaces.  We therefore begin by recalling some notions about b-pseudodifferential
operators and Sobolev spaces.

Let $\bPs m X$ denote the space of b-pseudodifferential (or ``totally
characteristic'') operators of order $m$ on $X$, developed by Melrose
\cite{MR83f:58073}.  This calculus of operators is designed to contain the
vector fields tangent to $\pa X$ as operators of order one.  The symbol of
an operator in the b-calculus naturally lies in the the b-cotangent bundle
$\Tbstar X$ introduced above.  Let $\bH m(X)$ denote the associated scale
of Sobolev spaces, defined as usual with respect to the measure $dx\,
dy/x$; for positive integral orders, these spaces may be defined by
regularity under application of products of vector fields of the form $x
D_x,$ $D_{y_i}.$

Let $(X,g)$ be a conic manifold.  Let $L^2_g (X)$ denote the space of
metric square-integrable functions on $X$ and let $\dcal_s= \Dom
(\Lap^{s/2})$, where $\Lap$ is the Friedrichs extension of the conic
Laplacian on $X$ acting on $L^2_g$.

The following is proven in \cite{Melrose-Wunsch1}:
\begin{proposition}
For $s \in (-n/2,n/2),$ $\dcal_{s}=x^{-n/2+s} \bH s (X)$.
\end{proposition}
\noindent (Note that the factor $x^{-n/2}$ arises because $L^2_g(X)= x^{-n/2}\bL(X).$)

The Cauchy problem for the wave equation is well-posed on Cauchy data in
$\dcal_s\oplus \dcal_{s-1}$ for any $s\in \RR$ and yields a solution $u\in
\mathcal{C}^0(\RR; \dcal_s) \cap \mathcal{C}^1(\RR; \dcal_{s-1}).$ It is
convenient to introduce the family of operators $\Theta_{r}$, $r \in
\RR\backslash \ZZ$ defined by
$$
\kappa(\Theta_r)(t,t') = \psi(t-t') \kappa(\abs{D_t}^r)(t,t')
$$
where $\psi(t)$ is a smooth function of compact support, equal to $1$ near
$t=0$, and $\kappa$ denotes Schwartz kernel.  (Hence $\Theta_r$ is a
version of $\abs{D_t}^r,$ but with Schwartz kernel cut off to have proper
support.)  Using the functional calculus, we easily see that $\Theta_r$
maps a solution to the conic wave equation $u\in \mathcal{C}^0(\RR;
\dcal_s) \cap \mathcal{C}^1(\RR; \dcal_{s-1})$ to a solution in
$\mathcal{C}^0(\RR; \dcal_{s-r}) \cap \mathcal{C}^1(\RR; \dcal_{s-r-1}).$
Moreover, $\Theta_r \Theta_{-r} u -u\in \mathcal{C} (\RR; \dcal_\infty).$
(See \cite{Melrose-Wunsch1} for details.)

Let $\trclass(H)$ denote the space of trace class operators on the Hilbert
space $H.$ The following lemma will be crucial in what follows, as it
yields a criterion for an operator to be trace-class in terms of its
mapping properties on domains.

\begin{lemma}\label{lemma:domaintrace}
If $A: \dcal_s \to \dcal_{s+r}$ with $r>n$ then $A \in \trclass(\dcal_s).$
\end{lemma}
\begin{proof}
To begin with we may reduce to the case in which $s=0,$ hence
$\dcal_s=L^2_g,$ by using the isomorphism $\ang{\Lap}^{s/2}: \dcal_{s'} \to
\dcal_{s'-s}.$

It follows from results of \cite{Melrose-Wunsch1} that $A: x^{n/2}\bL(X)
\to x^{-\ep} \bH {n+\ep} (X)$.  Let $B$ be an elliptic element of $\bPs
{n+\ep} X$, and let $C\in \bPs {-n-\ep} X$ be a ``small parametrix'' for
$B$ as constructed in \cite{Melrose:APS}, so that $CB -I=R \in \bPs
{-\infty} X.$ We now write
$$
A=(C x^{n/2-\ep}) (x^{-n/2+\ep} B A) + (Rx^\ep) (x^{-\ep}A).
$$
Proposition~4.57 of \cite{Melrose:APS} states that if $\alpha>0$ and
$\beta<-n$ then $x^\alpha \bPs \beta X \subset \trclass(L^2_g).$ Hence $C
x^{n/2-\ep}$ and $R x^\ep$ are of trace-class; since $x^{-n/2+\ep} BA$ and
$x^{-\ep} A$ are bounded on $L^2_g$, we conclude that $A$ is of trace
class.
\end{proof}

Finally, we recall from \cite{Melrose-Wunsch1} the propagation results
which, together with Lemma~\ref{lemma:domaintrace}, will enable us to
evaluate the regularity of the wave trace.
Let
$$
\prop(t)=
\begin{pmatrix}
\cos t\sqrt\Lap && \sin t \sqrt\Lap/\sqrt\Lap \\ - \sqrt\Lap \sin t
\sqrt\Lap && \cos t\sqrt\Lap
\end{pmatrix}
$$
be the solution operator to the Cauchy problem for the wave equation on
$X$.  Let $\Psi_c^* (X^\circ)$ denote the space of pseudodifferential
operators on $X^\circ$ with compact support.  Let $L(V_1, V_2)$ denote the
space of continuous linear maps between Fr\'echet spaces $V_i$.  Let
$$\E_s=\dcal_s\oplus \dcal_{s-1}$$ denote the ``energy space,'' with
$\E_\infty = \bigcap_{s\in\RR} \E_s$, $E_{-\infty}= \bigcup_{s\in\RR}\E_s.$

\begin{proposition}\label{prop:diff1}
Let $I\subset \RR_+$ be contained in a small neighborhood of $0.$ Let $A_1,A_2
\in \Psi^0(X^\circ)$ be supported near $\pa X$, and assume that that there
do not exists points $p\in \WF' A_2$, $q\in \WF' A_1$, and $t\in I$ such that
$p \dtilde q$.  Then $A_1 U(t) A_2 \in \mathcal{C}(I; L(\E_{-\infty}, \E_{\infty})).$
\end{proposition}

In the following proposition, $x$ and $x'$ will denote the boundary
defining function (in the product-type coordinates introduced in
\S\ref{section:geodesic}) in the left and right factors of the Schwartz
kernel.
\begin{proposition}\label{prop:geom}
Let $I$, $A_1$, $A_2$ be as in Proposition~\ref{prop:diff1}, and assume that
that there do not exists points $p\in \WF' A_2$, $q\in \WF' A_1$, and $t\in
I$ such that $p \gtilde q$.  Then for $t\in I$,
$$
A_1 U(t) A_2 = \begin{pmatrix} u_1 & u_2 \\ -u_2 & u_1 \end{pmatrix}
$$
where $u_1$ and $u_2$ are conormal distributions at $t=x+x'$ in
$H^{-1/2-\ep}$ resp.\ $H^{1/2-\ep}$ for all $\ep>0$.
\end{proposition}
A consequence of this proposition which we will use below is that subject
to the hypotheses, the diagonal term $u_1$ of $A_1U(t)A_2$ is a Fourier
integral operator of order $-\ep$ in the variable $x$, with values in
$\Psi^{-\infty}(Y)$, uniformly in $t$.  (Recall that the fundamental
solution discussed in the introduction is the off-diagonal term.)

We will also need a version of Proposition~\ref{prop:diff1} which applies
up to the boundary of $X$.
\begin{proposition}\label{prop:diff2}
Proposition~\ref{prop:diff1} still holds if either or both of $A_1$, $A_2$
are replaced by the operator of multiplication by a smooth function $\psi_i$
supported in a neighborhood of a component of $\pa X,$ and $\WF' A_i$
correspondingly replaced by $\pi^{-1}(\supp \psi_i).$
\end{proposition}

\section{Proof of the theorem}

As the proof has several parts, with ramifying subcases, this section is
similarly subdivided.  We begin with a sketch of the proof.  We then
construct the microlocal partition of unity whose properties are the
central ingredient in the proof.  Finally we prove successively the overall
regularity of the trace of the wave operator, the smoothness away from
lengths of diffractive geodesics, and the regularity away from lengths of
geometric geodesics.

The strategy of the proof is an imitation of the following non-constructive
proof of the Poisson relation for a compact boundaryless manifold $X$.  Let
$1=\sum A_i^* A_i+R$ be a microlocal partition of unity, with the $A_i$'s
having microsupport in very small sets in $S^*X$ and $R\in
\Psi^{-\infty}(X).$ Then if we let $\equiv$ denote equivalence modulo
$\CI(\RR)$, we have
\begin{equation}
\Tr D_t^{2k} U(t) \equiv \Tr \sum A_i^* A_i D_t^{2k} U(t) = \sum \Tr A_i D_t^{2k}
U(t) A_i^* = \sum \Tr A_i \Lap^k U(t) A_i^*.
\label{boundaryless}
\end{equation}
Let $I\subset \RR_+$ be a small interval containing no lengths of closed
geodesics.  Provided the microlocal partition of unity is sufficiently
fine, there are no bicharacteristics of length $t\in I$ with both endpoints
in $\WF' A_i.$ Hence by H\"ormander's propagation theorem, each term $A_i
\Lap^k U(t) A_i^*$ maps distributions to smooth functions, therefore is of
trace class.  Hence the wave trace is smooth on $I.$ Note that to make this
argument rigorous, we should actually approximate the $A_i$ by smoothing
operators in order to justify the formal manipulations in
\eqref{boundaryless}.

The simple approach sketched above requires modification on a conic
manifold.  First, owing to the global nature of the ``indicial operator''
at the boundary (one of the obstructions to compactness for operators in
Melrose's b-calculus), we must settle for a partition of unity of the form
$\sum_j \psi_j^2 +\sum_i A_i^* A_i$ where the $\psi_j$'s are smooth cutoff
functions localizing at different boundary components and the $A_i$'s have
microsupport over a relatively compact subset of $X^\circ.$ Furthermore, we
only have a solid understanding of propagation of singularities for a short
time, in which singularities cannot repeatedly interact with $\pa X.$ Hence
we will decompose $t=t_0+\dots t_L$ with each term small, and write
\begin{multline*}
\Tr U(t)\equiv \\\Tr (\sum_j \psi_j^2 +\sum_i A_i^* A_i) U(t_0) (\sum_j \psi_j^2 +\sum_i
A_i^* A_i) U(t_1) \dots (\sum_j \psi_j^2 +\sum_i A_i^* A_i) U(t_L).
\end{multline*}
We rewrite this sum as a single sum over terms of the form
$$
\Tr (B_{i_0}^* B_{i_0}) U(t_0)(B_{i_1}^* B_{i_1}) U(t_1)  \dots (B_{i_L}^* B_{i_L}) U(t_L)
$$
where each $B_i$ is either one of the $\psi_j$'s or one of the $A_i$'s.  The
strategy for the part of the theorem dealing with diffractive geodesics
is to show that if the partition of unity is taken sufficiently fine, one
factor in each summand of the form $B_{i_j} U(t_j) B_{i_{j+1}}^*$ (possibly
after cyclic permutation), is smoothing.  This is based on the absence of
closed diffractive geodesics of length $t$, and on the principle that if
there are no closed geodesics of length $t$ then there are no closed
``broken'' geodesics of such length which are repeatedly allowed to
propagate for time $t_j$ and then jump arbitrarily within $\WF'B_{i_j}.$

The fact that our partition of unity fails to localize over the boundary is
of no import in dealing with diffractive geodesics which may, after all,
jump freely from one point over the boundary to another in the same
component.  Difficulties arise, however, in dealing with geometric
geodesics.  In particular, the strategy outlined above does not work in the
geometric case, and instead of considering terms of the form $B_{i_j}
U(t_j) B_{i_{j+1}}^*,$ if one of the $B$'s is a localizer near a boundary
component we may need to consider more complicated factors $B_{i_j} U(t_j)
B_{i_{j+1}}^*B_{i_{j+1}} U(t_{j+1}) B_{i_{j+2}}^*$ so as to be considering
interior-to-interior propagation of singularities.  This situation requires
the consideration of a number of cases, involving various possibilities for
existence of geodesics from $\WF' B_{i_{j+2}}$ to $\WF' B_{i_{j+1}}$ and
thence to $\WF' B_{i_{j}}.$

\subsection{Microlocal partition of unity}
The reader may find it helpful to refer back and forth between this section
and those following in which its results are employed to prove the theorem.

For $j=1,\dots, K$, let $U_j$ denote a product neigborhood of $Y_j$ of the
form $\{x<\ep_j\}$; let $\psi_j\in \CI(X)$ denote a cutoff function
supported in $U_j$, equal to $1$ on a smaller neighborhood
$U_j'=\{x<\ep_j'\}$ of $Y_j$.  Let $\{V_i\}_{i=1}^M$ be an open cover of $S^*
X\backslash \pi^{-1} \bigcup_i U_i'$.

We now establish the existence of operators $A_i \in \Psi_c^0(V_i)$
such
that
$$
R+\sum_{j=1}^K \psi_j^2+ \sum_{i=1}^M A_i^* A_i =1.
$$
with $R \in \Psi_c^{-\infty}(X^\circ)$ a compactly supported smoothing
operator.  ($\Psi_c^0(V_i)$ denotes the space of pseudodifferential
operators whose kernels have supports which project to compact subsets of
$V_i$ on both left and right factors.)  To do this note that using a
symbolic construction in the ordinary pseudodifferential calculus, we can
choose $D_i \in \Psi_c^0(V_i)$ such that $\WF' \sum D_i^* D_i -1 \subset
\bigcup U_j'$.  Then setting $\zeta=(1-\sum \psi_j^2)^{1/2}$, we may set
$A_i=D_i \zeta$, hence $\sum A_i^* A_i -\zeta^2\in \Psi^{-\infty}_c
(X^\circ),$ as desired.

Let $A_{i,\delta}$, $\delta\in [0,1]$, be a family in
$\Psi^{-\infty}(X^\circ)$ converging strongly to $A_i$ as $\delta
\downarrow 0$.

Let 
$$
B_{i,\delta}=
\begin{cases}  \psi_i \text{ for } i=1,\dots, K,\\ A_{i-k,\delta}
  \text{ for } i=K+1,\dots K+M;
\end{cases}
$$
 let $B_i=B_{i,0}$. Hence
$$
1=\sum B_i^* B_i + R.
$$

Let $\bar t\in \RR_+$ be sufficiently small that $\{d(Y_i, \cdot)<\bar t\}$
is contained in a product neighborhood of $Y_i$ for all $i$, and such that
$\bar t \ll d(Y_i, Y_j)$ for all $i,j$.

Suppose that we are given $T \in \RR_+\backslash \dif.$ Decompose
$T=t_0+\dots+t_L$ such that $\bar t /4 <t_i<\bar t/2$ for all $i$ (we may
decrease $\bar t$ if necessary).

\emph{In the statement of the following lemma and thenceforth, we always
consider indices $i_0, \dots, i_L$ up to cyclic shift, i.e.\ $i_l$ is
always shorthand for $i_{l \bmod{L+1}}$.}
\begin{lemma}
Given $T=t_0+\dots t_L \in \RR_+\backslash \dif$ with each $t_i \in (\bar
t/4, \bar t/2),$ we may choose the cover $\cvr=\{\pi^{-1} U_i\}\cup
\{V_i\}$ of $\Sbstar X$ sufficiently fine that for any ${i_0},\dots, {i_L}$
there exist $l$ and an open interval $I \ni t_{i_l}$ such that there does
not exist a lifted diffractive geodesic with length in $I$, ending in
$\WF'B_{i_l},$ and beginning in $\WF' B_{i_{l+1}}.$
\label{lemma:nondiffractive}
\end{lemma}

\begin{proof}
Let $\cvr^m$ denote a sequence of finer covers $\{\pi^{-1} U_i^m\}\cup
\{V_i^m\},$ and $B_l^m$ the corresponding operators in our microlocal
partition of unity.  If the conclusion of the lemma fails to holds then for
each $m$ there exists an $(L+1)$-tuple of pairs of points $(p_0^m, q_0^m),
\dots, (p_L^m, q_L^m)$ such that $p_l^m, q_l^m \in \WF' B_{i_l}$ and $t_l^m
\to t_l$ such that $p_{l+1}^m \dtilde[t_l^m] q_l^m.$ As $m\to \infty$,
i.e.\ as the partition is refined, we can then extract a subsequence such
that $p_l^m \to \bar p_l$ and $q_l^m \to \bar q_l.$ If either $\bar p_l$ or
$\bar q_l$ is in $\Sbstar X^\circ$, then $\bar p_l=\bar q_l$, since the two
points lie in a sequence of shrinking sets over $X^\circ$.  If $\bar p_l$
or $\bar q_l$ is in $\Sbstar_{\pa X} X$ then the two points $\bar p_l$ and
$\bar q_l$ must lie over the same component of $\pa X$.  Thus by
Proposition~\ref{prop:closed}, each $\bar p_{l+1}$ is connected to $\bar
p_l$ by a diffractive geodesic of length $t_l$, hence these points must
lie along a closed diffractive geodesic of length $T$ (recall that $l$ is
considered modulo $L$ throughout), contradicting our assumptions.
\end{proof}

In the case $T \in \RR_+ \backslash \geo$, we will need to impose a subtler
set of conditions on our refined cover.
\begin{lemma}
Given $T=t_0+\dots t_L \in \RR_+\backslash \geo$ with each $t_i \in (\bar
t/4, \bar t/2),$ we may choose the cover $\{\pi^{-1} U_i\}\cup \{V_i\}$ of
$\Sbstar X$ sufficiently fine that for any $i_0,\dots, i_L$ there exist
open intervals $I_l \ni t_l$ such one of the following holds:
\begin{enumerate}
\item there exists $l$ such that $i_l$ and $i_{l+1}$ are both less than or
equal to $K$ (i.e.\ are both $B_{i_l}$ and $B_{i_{l+1}}$ are $\psi$'s), or
\item there exists $l$ such that $i_l$ and $i_{l+1}$ are both greater than
$K$ and there does not exist a geometric geodesic with length
in $I_l$, ending in $\WF'A_{i_l-K},$ and beginning in $\WF'
A_{i_{l+1}-K},$ or
\item there exists $l$ such that $i_l,\ i_{l+2}>K$ and $i_{l+1}\leq K$ and
there does not exist a geometric geodesic with length in $I_l +
I_{l+1},$ ending in $\WF' A_{i_l-K}$ and beginning in $\WF' A_{i_{l+2}-K}.$
\end{enumerate}
\label{lemma:nongeometric}
\end{lemma}

\begin{proof}
If the conclusion fails to hold, there must exist a sequence of shrinking
covers $\cvr^m$ and shrinking intervals $I_l^m \ni t_l$ in each of
which there exists a ``word'' $i_0^m, \dots i_L^m$ and $(p^m_l, q^m_l) \in
\WF' B_{i_l}^m$ such that such that no two successive
$i_l$'s are less than or equal to $K$, and for all $l=1, \dots L$, either
\begin{itemize}
\item
if $i^m_l, i^m_{l+l} >K$, there exists a lifted geometric geodesic from
$p_{i_{l+1}}^m$ to $q_{i_l}^m$ with length in $I_l^m$
or
\item
if $i^m_l> K$ and $i^m_{l+1}\leq K,$ then $i^m_{l+2}>K$ and there exists a
lifted geometric geodesic from $p_{i_{l+2}}^m$ to $q_{i_l}^m$ with length
in $I_l^m+I_{l+1}^m.$
\end{itemize}

As $m\to \infty$, i.e.\ as the partition is refined, we may extract a
subsequence such that $p_l^m \to \bar p_l$ and $q_l^m \to \bar q_l$.  If
eventually $i_l^m>K$, i.e.\ if $p_l^m$ and $q_l^m$ eventually lie in one of
the interior sets $V_k^m$, then $\bar p_l=\bar q_l$.  Hence if both $p_l^m$ and
$p_{l+1}^m$ are such points then $\bar p_{l+1} \gtilde[t_l] \bar p_l$ by
Proposition~\ref{prop:closed}.  If, on the other hand, we eventually have
$i_l^m>K$ but $i_{l+1}^m\leq K$, then $i_{l+2}^m >K,$ and we must have
$\bar p_{l+2}\gtilde[t_l+t_{l+1}] \bar p_l$ by
Proposition~\ref{prop:closed}.  Thus, by Proposition~\ref{prop:closed}, the
points $\bar p_l$ must all lie along a geometric geodesic of length $T$,
contradicting $T \notin \geo.$
\end{proof}

\begin{lemma}
By further refining the cover $V_j$ and shrinking the intervals $I_j$ we
may assume in case 3 of Lemma~\ref{lemma:nongeometric} that either
\begin{enumerate}\renewcommand{\theenumi}{\roman{enumi}}
\item
there exist no \emph{diffractive} geodesics with length in $I_l$ ending in
$\WF'A_{i_l-K},$ and beginning in $\pi^{-1} \supp \psi_{i_{l+1}}$ or 
\item
there exist no \emph{diffractive} geodesics with length in $I_{l+1}$ ending in
 $\pi^{-1} \supp \psi_{i_{l+1}}$ and beginning in $\WF'A_{i_{l+2}-K}$ or
\item
any diffractive geodesic of length $\mathcal{L} \in I_l+I_{l+1}$ beginning in
$\WF'A_{i_{l+2}-K}$ and ending in $\WF'A_{i_{l}-K}$ must lie in $\pi^{-1}
\supp \psi_{i_{l+1}}$ either for all $t \in I_{l+1}$ or for all $t \in \mathcal{L}-I_l$
or
\item
no diffractive geodesics with length in $I_l$ ending in
$\WF'A_{i_l-K}$ pass through $\pa X$
or
\item
no diffractive geodesics with length in $I_{l+1}$ beginning in
$\WF'A_{i_{l+2}-K}$ pass through $\pa X$.
\end{enumerate}
\label{lemma:subcases}
\end{lemma}
\begin{proof}
Since $\bar t \ll \min d(Y_i, Y_j)$ for any given values of $i_l$,
$i_{l+2}$ there is at most one value of $i_{l+1}$---call it $\bar \imath$---for
which (i) and (ii) fail to hold (provided the intervals $I_l$ are chosen
sufficiently small).  If $I_l$ is a sufficiently small interval, the set
$\bigcap_{t \in \bar I_l} \{p: q\dtilde[t] p \Rightarrow q \in \pi^{-1}
U_{\bar \imath}' \}$ is an open neighborhood of the compact set $\bigcup_{t \in
\bar I_l} \{p : q \dtilde[t] p \text{ for some } q \in \pi^{-1} Y_{\bar
i}\}$ (cf.\ Lemma~4.5 of \cite{Wunsch2}), hence we may refine the cover
$\{V_j\}$ so that any $V_j$ intersecting the latter lies inside the former.
Similarly, we may arrange that any $V_j$ intersecting $\bigcup_{t\in \bar
I_l} \{p : p \dtilde[t] q \text{ for some }q \in \pi^{-1} Y_{\bar \imath}\}$
lies inside $\bigcap_{t \in \bar I_l} \{p: p \dtilde[t] q \Rightarrow q \in
\pi^{-1} U_{\bar \imath}' \}.$

Any diffractive geodesic with length in $I_l+ I_{l+1}$ ending in
$\WF'A_{i_l-K},$ and beginning in $\WF' A_{i_{l+2}-K}$ must pass through
the boundary component over $Y_{\bar \imath},$ as there are no geometric
geodesics connecting such points.  Let $\gamma$ be such a geodesic
(time-parametrized) and $\mathcal{L}$ its length.  If $\gamma(t)$ passes
through $Y_{\bar \imath}$ for some $t \in I_{l+1} \cup (\mathcal{L}-I_l)$ then
by our refinement of the cover in the previous paragraph, $\gamma(t) \in
\pi^{-1} U_{\bar \imath}'$ for all $t$ in either $I_{l+1}$ or $\mathcal{L}-I_l$,
hence case (iii) holds.

If, on the contrary, every geodesic $\gamma(t)$ with length in $I_l+
I_{l+1}$ ending in $\WF'A_{i_l-K},$ and beginning in $\WF' A_{i_{l+2}-K}$
passes through $Y_{\bar \imath}$ at time $\tau \notin I_{l+1} \cup
(\mathcal{L}-I_l)$, we have two cases, $\tau \gtrless I_{l+1} \cup
(\mathcal{L}-I_l).$ We may assume, by further refinement of the cover, that
just one of these inequalities holds for all the aforementioned diffractive
geodesics with a given $i_l,\ i_{l+2}$; indeed we may assume that the set
of lengths of geodesics connecting a $V_j$ to a boundary component is a
single interval.  Say we are in the case $>$.  Then all geodesics from
$\WF' A_{i_{l+2}-K}$ to $\pi^{-1}Y_{\bar \imath}$ have length greater than
$I_{l+1}$, i.e.\ case (v) holds.  The proof with sign $<$ is analogous:
case (iv) holds.
\end{proof}

We will need one further piece of geometric information about the partition
of unity.
\begin{lemma}
Let $\bar t$ be fixed as above.
If the sets $U_j \supset Y_j$ are chosen sufficiently small, then there are
no diffractive geodesics of length in $[\bar t /4, \bar t /2]$
connecting points in $U_j$ with points in $U_k$ for any $j,k$.
\label{lemma:noloops}
\end{lemma}
\begin{proof}
The result is clear for $j \neq k$, since $\bar t$ is less than the
distance between boundary components.  Thus we treat the case $j=k$.  It
suffices to show that for any $j$, there exists $\ep$ sufficiently small
that any diffractive geodesic starting in the component of $\{x<\ep\}$
containing $Y_j$ must lie in $\{x>\ep\}$ for $t \in [\bar t/4, \bar t/2]$.
Any diffractive geodesic passing through $Y_j$ is locally of the form
$y=y_0,\ x=x_0\pm t$; along such a diffractive geodesic, $t=x\pm x_0$,
hence either $x$ or $x_0$ must exceed $\bar t/8$ when $t>\bar t/4$.  Hence
it suffices, in considering such geodesics, merely to take $\ep=\bar t/8.$

We now consider geodesics which do not pass through $Y_j$.  We consider the
unit speed geodesic flow to be the flow given by \eqref{spray} in $\Tbstar
X$, inside $g=(\xi^2+h(x,y,\eta))/x^2=1$.  If $x_0<\ep$, then of course
$\abs{\xi_0}<\ep$.  On the other hand $\dot \xi =(\xi^2+h)/x^2 +
\ocal(\eta^2/x) \geq 1/2$ if $\ep$ is sufficiently small.
Hence $\xi(t) \geq \xi_0+t/2 \geq -\ep+ t/2$.  Since $\dot x =\xi/ x$, we
obtain $x^2 \geq t^2/2-2\ep t$, which is greater than $\ep^2$ for $t>\bar
t/4$ if $\ep$ is sufficiently small.
\end{proof}

\subsection{Overall regularity}
We begin the proof of the theorem by establishing the overall regularity
$$
\Tr \prop(t) \in \cC^{-n-0}(\RR).
$$
To do this, note that for $\ep>0$ and all $s \in \RR,$
$$
\Theta_{-n-\ep} \prop (t): \E_s \to \E_{s+n+\ep},
$$
hence by Lemma~\ref{lemma:domaintrace}, $\Theta_{-n-\ep} \prop (t)$ is
locally in $L^\infty(\RR; \trclass (\E_s))$, i.e.\ $\Tr \prop(t) \in
\cC^{-n-0}(\RR)$ as desired.

\subsection{Diffractive lengths}

We now show that $\prop(t) \in \CI(\RR\backslash \dif)$.  It suffices by
selfadjointess of $\Lap$ to show smoothness in $\RR_+\backslash \dif.$

We now choose the cover $\cvr$ as guaranteed by
Lemma~\ref{lemma:nondiffractive}, and compute
\begin{multline*}
\Tr D_t^{k(L+1)} \prop (t+t_1+\dots+t_n) = \lim_{\delta\downarrow 0} \Tr
\big ( R+ \sum B_{i,\delta}^* B_{i,\delta}\big ) D_t^k \prop(t) \\ \cdot \big ( R+
\sum B_{i,\delta}^* B_{i,\delta}\big ) D_t^k \prop(t_1) \big ( R+ \sum
B_{i,\delta}^* B_{i,\delta}\big ) \dots  \big ( R+ \sum B_{i,\delta}^* B_{i,\delta}\big )
D_t^k \prop(t_L)
\end{multline*}
in the sense of distributions (since $1-\sum B_{i,\delta}^2-R$ approaches
zero in norm as a map $\dcal_s\to \dcal_{s+\ep}$ for all $\ep>0$).  Now
expand out the above sum.   By Lemma~\ref{lemma:nondiffractive}, each
term of the form
$$
B_{i_0,\delta}^* B_{i_0,\delta} D_t^k \prop(t) B_{i_1,\delta}^*
B_{i_1,\delta} \dots B_{i_L,\delta}^* B_{i_L,\delta} D_t^k \prop(t_L)
$$
contains (after a possible cyclic permutation) a factor $B_{i_l,\delta}
D_t^k \prop (t_{i_l}) B_{i_{l+1},\delta}^*$ that, by
Propositions~\ref{prop:diff1} and \ref{prop:diff2}, maps $\E_{-\infty} \to
\E_{\infty}$ (uniformly in $\delta$), while the rest
of the terms map $\E_s \to \E_{s-k}$ for all $s\in \RR$. Hence the
expression as a whole is of trace class, uniformly in $t$ near $t_0$ and as
$\delta \downarrow 0$.  Any term involving a factor of $R$ also has this
property, since $R: \E_{-\infty} \to \E_\infty.$ Thus, for all $k \in
\NN_0$, $D_t^{k(L+1)} \prop(t+t_1+\dots+t_L) \in L^\infty (\RR;
\trclass(\E_s))$ near $t=t_0$, proving that
$$
\prop(t) \in \CI(\RR\backslash \dif).
$$

\subsection{Geometric lengths}

It now remains to show
$$
\prop(t) \in \cC^{-1-0}(\RR\backslash \geo).
$$
Let $T \in \RR_+ \backslash \geo$ and decompose $T$ into a sum as before.
We have for $\ep>0$,
\begin{equation}
\begin{aligned}
\Theta_{-1-\ep} \prop (t+t_1+\dots+t_L) &= \Theta_{-1-\ep} \prop(t)\prop(t_1) \dots \prop(t_L)\\
&=
\big( R+ \sum B_{i,\delta}^* B_{i,\delta} \big)\Theta_{-1-\ep} \prop(t)\big(
R+ \sum B_{i,\delta}^* B_{i,\delta} \big)\prop(t_1)\\ & \dots \big( R+\sum B_{i,\delta}^* B_{i,\delta} \big)\prop(t_L),
\end{aligned}
\label{expansion}
\end{equation}

For brevity's sake, let
$$
W_j =
\begin{cases}
U(t_j), & j \neq 0\\
\Theta_{-1-\ep} U(t), & j=0.
\end{cases}
$$
Applying Lemma~\ref{lemma:nongeometric}, we may refine the cover of
$\Sbstar X$ sufficiently that any word
$$
B_{i_0}^* B_{i_0} W_0  B_{i_1}^* B_{i_1} W_1
B_{i_2}^* B_{i_2} \dots  B_{i_L}^* B_{i_L} W_L,
$$
in \eqref{expansion} contains (modulo cyclic permutation) a factor of one
of the following forms, corresponding to the various cases of the lemma:
\begin{enumerate}
\item 
$\psi_{i_l} W_l \psi_{i_{l+1}}$, or
\item 
$A_{i_l-K} W_l A_{i_{l+1}-K}^*$
where there does not exist a geometric geodesic with length
in $I_l$ ending in $\WF'A_{i_l-K},$ and beginning in $\WF'
A_{i_{l+1}-K},$ or
\item 
$A_{i_l-K} W_l \psi_{i_{l+1}}^2 W_{l+1} A_{i_{l+2}-K}^*$ where there does not
exist a geometric geodesic with length in $I_l+ I_{l+1}$ ending in
$\WF'A_{i_l-K},$ and beginning in $\WF' A_{i_{l+2}-K};$
\end{enumerate}

The term in case 1 above is certainly of trace class, since by
Lemma~\ref{lemma:noloops} there are no diffractive geodesics (let alone
geometric ones) connecting points in $\psi_{i_l}$ and $\psi_{i_{l+1}}$ and
having length in $[\bar t/4, \bar t/2]$.  Hence any term
$$
B_{i_0,\delta}^* B_{i_0,\delta} \Theta_{-1-\ep} \prop(t) B_{i_1,\delta}^* B_{i_1,\delta} \prop(t_1)
B_{i_2,\delta}^* B_{i_2,\delta} \dots B_{i_L,\delta}^* B_{i_L,\delta} \prop(t_L),
$$
containing such a factor is uniformly of trace class in $t$ near $t_0$ and
as $\delta \downarrow 0$.

Proposition~\ref{prop:geom} shows that a term of the type in case
2,
$$
A_{i_l-K} U(s) A_{i_{l+1}-K}^*,
$$
where there are no geometric geodesics with lengths in $I_l$ between the
microsupports of the outer terms, is in fact a system of Fourier integral
operators in $x$ of order $+\ep$ with values in $\Psi^{-\infty}(Y_j)$ for
some $j$, uniformly in $s \in I_l.$ (We are assuming for the moment that $l
\neq 0$.)  Now let $D \in \Psi^{-1-\ep}(X^\circ)$ be compactly supported
and elliptic on $\WF' A_{i_l-K};$ let $E \in \Psi^{+1+\ep}(X^\circ)$ be a
compactly supported microlocal parametrix for $D$, hence $\WF' E D \cap
\WF' A_{i_l-K}=\emptyset$.  Then
$$
A_{i_l-K} U(t_l)A_{i_{l+1}-K}^* = E D A_{i_l-K} U(t_l)
A_{i_{l+1}-K}^* + S
$$
where $S$ is a compactly supported smoothing operator on $X^\circ$, hence
does not contribute to singularities of the trace.  The
term
$$
D A_{i_l-K,\delta} U(t_l)
A_{i_{l+1}-K,\delta}^*
$$
is then of trace class, uniformly in $\delta$; hence the whole term
$$
\Tr B_{i_0,\delta}^* B_{i_0,\delta} \Theta_{-1-\ep} \prop(t)
B_{i_1,\delta}^* B_{i_1,\delta} \prop(t_1) 
B_{i_2,\delta}^* B_{i_2,\delta} \dots  B_{i_L,\delta}^* B_{i_L,\delta} \prop(t_L)
$$
can be written, modulo an error term mapping $\E_{-\infty} \to
\E_\infty,$ as a product of factors preserving $\E_s$ for all $s$,
times a factor $\Theta_{-1-\ep} U(t)$ which maps $\E_s \to
\E_{s+1+\ep}$ for all $s\in \RR$ and $\ep>0$, times a factor $E$ which
maps $\E_s \to \E_{s-1-\ep}$, times a factor
$$
D A_{i_l-K,\delta} U(t_l)
A_{i_{l+1}-K,\delta}^* \in \trclass(\E_s),
$$
uniformly in $\delta$.  Hence this term has a trace uniformly bounded as
$\delta\downarrow 0$.  The terms involving the smoothing error $R$ do not
contribute, as before.  The separate case $l=0$ works similarly.

In the subcases (i) and (ii) of case (3) (described in
Lemma~\ref{lemma:subcases}) the corresponding term in the decomposition of
$U(t)$ again maps $\E_{-\infty} \to \E_\infty$, hence does not
contribute to singularities of the trace.  To deal with subcase (iii), we
decompose
\begin{multline*}
A_{i_l-K} U(t_l) \psi_{i_{l+1}}^2 U(t_{l+1}) A_{i_{l+2}-K}^* \\ =A_{i_l-K} U(t_l+t_{l+1}) A_{i_{l+2}-K}^*
+ A_{i_l-K} U(t_l) (1-\psi_{i_{l+1}}^2) U(t_{l+1}) A_{i_{l+2}-K}^*.
\end{multline*}
The second term maps $\E_{-\infty}\to \E_\infty$, while the first can
be dealt with just as the similar term in case (2).  For case (iv), note
that by Egorov's theorem, there exists $\tilde A \in \Psi_c^{0} (X^\circ)$
with microsupport in the geodesic flowout for time $-t_l$ of $\WF'
A_{i_l-K}$ such that
$$
A_{i_l-K} U(t_l) \psi_{i_{l+1}}^2 U(t_{l+1}) A_{i_{l+2}-K}^* = U(t_l)
\tilde A U(t_{l+1}) A_{i_{l+2}-K}^*.
$$
Now by assumption, there are no geometric geodesics of length in $I_{l+1}$
beginning in $\WF' A_{i_{l+2}-K}$ and ending in $\WF' \tilde A$, hence we
may proceed as in case (2), since by \cite{Melrose-Wunsch1}, the term
$\tilde A U(t_{l+1}) A_{i_{l+2}-K}^*$ is a Fourier integral operator in $x$
of order $+\ep$ with values in $\Psi^{-\infty}(Y_j)$ for some $j$.  Case
(v) of Lemma~\ref{lemma:subcases} is analogous.

\bibliography{All}
\bibliographystyle{amsplain}
\end{document}